\documentclass{amsart}
\usepackage{amssymb,latexsym}
\usepackage{amsmath}
\usepackage{graphicx}

\newtheorem{thm}{Theorem}[section]
\newtheorem{cor}[thm]{Corollary}
\newtheorem*{thm1.4}{Theorem 1.4}

\newtheorem{defin}[thm]{Definition}
\newtheorem{lemma}[thm]{Lemma}

\newtheorem{prop}[thm]{Proposition}
\newtheorem{rmk}[thm]{Remark}

\newcommand{\bdd}{\mbox{$\partial$}}

\newcommand{\mcA}{\mbox{$\mathcal{A}$}}
\newcommand{\mcB}{\mbox{$\mathcal{B}$}}
\newcommand{\mcU}{\mbox{$\mathcal{U}$}}
\newcommand{\mcV}{\mbox{$\mathcal{V}$}}
\newcommand{\Sss}{\mbox{$\Sigma$}}
\newcommand{\aaa}{\mbox{$\alpha$}}
\newcommand{\Ggg}{\mbox{$\Gamma$}}

\newcommand{\xX}{\mbox {\sc x}}
\newcommand{\yY}{\mbox {\sc y}}
\newcommand{\aA}{\mbox {\sc a}}
\newcommand{\bB}{\mbox {\sc b}}
\begin{document}

\title{Distance of Heegaard splittings of knot complements}

\author{Maggy Tomova}
\address{\hskip-\parindent
 Maggy Tomova\\
  Mathematics Department \\
Rice University \\
 Houston, TX 77005, USA}
\email{maggy.tomova@rice.edu}

\begin{abstract}
 Let $K$ be a knot in a closed orientable irreducible 3--manifold
 $M$ and let $P$ be a Heegaard splitting
 of the knot complement of genus at least two. Suppose $Q$ is a bridge surface for $K$ and let $N(K)$ denote a regular neighborhood of $K$. Then either
 \begin{itemize}
 \item $d(P)\leq 2-\chi(Q-(N(K))$, or
 \item $K$ can be isotoped to be disjoint from $Q$ so that after the isotopy $Q$ is
 a Heegaard surface for $M-N(K)$ that is isotopic to a possibly
 stabilized copy of $P$.

 \end{itemize}
\end{abstract}

\maketitle
\section{Introduction}

A Heegaard splitting of a compact 3--manifold $M$ is a
decomposition of the manifold into two compression bodies, $A$ and
$B$.  If the manifold is closed, $A$ and $B$ are handlebodies.
The common boundary of $A$ and $B$ is called a {\em Heegaard
surface} which we denote by $P$. We will write $M=A \cup_P B$.

The distance between any two essential simple closed curves
$\alpha$ and $\beta$ in a Heegaard surface $P$ is the smallest
integer $n\geq 0$ so there is a sequence of essential simple
closed curves $\alpha=\alpha_0,...,\alpha_n=\beta$ in $P$ such
that for each $1\leq i \leq n$, $\alpha_{i-1}$ and $\alpha_i$ can
be isotoped to be disjoint in $P$. The {\em distance} of a
Heegaard splitting $A \cup_P B$, $d(P)$, defined by Hempel
\cite{Hem}, is the smallest integer $n$ so that there is an
essential curve $\alpha$ in $P$ which bounds a disk in $A$ and an
essential curve $\beta$ which bounds a disk in $B$ and $d(\alpha,
\beta)=n$.  If $d(P)=0$ then we say that the Heegaard splitting $A
\cup_P B$ is {\em reducible}, if $d(P)=1$ then we say that $A
\cup_P B$ is {\em weakly reducible} and if $d(P)\geq 2$ then we
say that $A \cup_P B$ is {\em strongly irreducible}.

 Obtaining bounds on the distance of a Heegaard splitting $M=A \cup_P B$ is an interesting problem and much
 progress has been made recently.
  The first such bound is due to
 Hartshorn \cite{Hart} who uses a closed essential surface in the manifold to give a
 bound on $d(P)$. This result was extended by Scharlemann in \cite{Sc2}
 to allow for possible boundary components of the essential surface.

 \begin{thm}\label{thm:essentialsurface} \cite[Theorem 3.1]{Sc2}
Suppose $P$ is a Heegaard surface for a compact orientable
manifold $M$ and $(Q, \bdd Q)\subset (M,\bdd M)$ is a connected
essential surface. Then $d(P)\leq 2-\chi(Q)$

\end{thm}

 Scharlemann and Tomova, \cite{ST}, gave
a bound on $d(P)$ using a second Heegaard surface.

\begin{thm}\label{thm:alternate}\cite[Corollary 3.5]{ST}
Suppose $P$ and $Q$ are distinct Heegaard splittings for the
compact orientable $3$-manifold $M$.  Then either $d(P) \leq 2
-\chi(Q)$ or $Q$ is isotopic to a stabilization of $P$.
\end{thm}

Suppose a manifold $M$ contains a knot $K$ and let $M_K$ be the
knot complement. There has been interest in obtaining a bound on
the distance of Heegaard splittings of $M_K$ using a bridge
surface for the knot. We say that a surface $Q$ is a bridge
surface for a knot $K$ in a manifold $M$ if $Q$ is a Heegaard
surface for $M=X \cup_Q Y$ and $K$ intersects each of the
compression bodies $X$ and $Y$ in arcs that are simultaneously
parallel to $Q$.

Recently Johnson and Thompson \cite[Theorem 1]{JT} used Theorem
\ref{thm:alternate} to show that if a knot $K \subset M$ with
bridge surface $Q$ is such that $K$ cannot be isotoped to lie in
$Q$, then the distance of any Heegaard splitting of the knot
complement is bounded by $\chi(Q-N(K))$ where $N(K)$ is a regular
neighborhood of $K$. They used this theorem to obtain tunnel
number one knots that have arbitrarily high bridge number with
respect to any
 bridge surfaces of genus 1. Kobayashi and
 Rieck \cite[Proposition 2.6]{KR} extended this result to obtain the same bound on the
 distance of a Heegaard splitting in the case when the genus of the bridge surface $Q$ is strictly less
 than the minimum genus of any Heegaard splitting of the knot complement. Their
 result allows for the knot to be isotopic to an embedded curve in
 $Q$. They used this distance bound to prove the existence of counterexamples
 to Morimoto's conjecture \cite[Conjecture 1.5]{Mor}.

 In this paper we remove the restriction on the genus of $Q$ as well as
 restrictions on whether the knot is isotopic into $Q$. More precisely:

 \begin{thm} \label{thm:main}
 Let $K$ be a knot in a closed oriented irreducible 3--manifold
 $M$ and let $P$ be a Heegaard splitting
 of the knot complement of genus at least two. Suppose $Q$ is a bridge surface for $K$
 and let $N(K)$ denote a regular neighborhood of $K$. Then either
 \begin{itemize}
 \item $d(P)\leq 2-\chi(Q-N(K))$, or
 \item $K$ can be isotoped to be disjoint from $Q$ so that after the isotopy $Q$ is
 a Heegaard surface for the knot complement that is isotopic to a possibly
 stabilized copy of $P$.

 \end{itemize}
 \end{thm}

In some sense this theorem completes the solution to the problem
as it is clear that a stabilized copy of $P$ will not carry any
information about $d(P)$ and these are the only bridge surfaces
which we exclude. This result was used by Minsky, Moriah and
Schleimer in \cite{MMS} to show that for any pair of integers $t$
and $b$ there is a knot with tunnel number $t$ such that every
bridge surface of genus $t$ for the knot intersects the knot in at
least $2b+2$ points. The proof of Theorem \ref{thm:main} uses
different methods than those used in \cite{JT} and \cite{KR} and
is independent of them.

\section{Definitions and Notation} \label{sec:defin}

In this paper, unless otherwise specified, we will consider a
closed irreducible orientable 3--manifold $M$ containing a knot
$K$. If $X$ is any subset of $M$ we will denote by $X_K$ the set
$X-N(K)$ where $N(K)$ is an open tubular neighborhood of $K$. In
particular the knot exterior $M-N(K)$ will be denoted by $M_K$.
Note that $M_K$ is a compact orientable manifold with a single
torus boundary component. Throughout this paper $P$ will be a
surface splitting $M$ into handlebodies $A$ and $B$ so that $K$ is
entirely contained in $A$. We will further assume that $A_K$ is a
compression body, that is $P$ is also a Heegaard surface for $M_K$
and $M_K=A_K \cup_P B$. We will further assume that $P$ has genus
at least 2 as otherwise distance is always infinite with our
definition.

Recall that a simple closed curve in a compact surface $S$ is {\em
essential} if it does not bound a disk in the surface and it is
not parallel to a boundary component of the surface.  Suppose $M$
is a closed manifold containing a knot $K$ and $S$ is a surface in
$M$ transverse to $K$. A disk $D \subset M_K$ is a compressing
disk for $S_K$ if $D\cap S=\bdd D$ and $\bdd D$ is an essential
simple closed curve in $S_K$. A disk $D^c$ in $M$ is a {\em cut
disk} for $S_K$ if $D^c \cap S_K=\bdd D^c$, $\bdd D^c$ is
essential in $S_K$ and $D^c$ intersects $K$ in a single point. A
{\em c-disk} is a cut or a compressing disk.

    A properly embedded surface $F$ in a 3--manifold is {\em essential}
    if it is incompressible and has at least one component that is not
    parallel to $\bdd M$.

A surface $F$ in $M$ is called a {\em splitting surface} if $M$
can be written as the union of two 3--manifolds $U$ and $V$ along
$F$. If $F$ is a splitting surface it is {\it bicompressible} if
it is compressible in both $U$ and $V$ and $F$ is {\it
c-bicompressible} if it has c-disks in both $U$ and $V$. If $F$ is
a splitting surface for $M$, we will call $F$ {\em c-weakly
incompressible} if it is c-bicompressible and any pair of c-disks
for $F$ on opposite sides of the surface intersect along their
boundaries. If a c-bicomressible surface $F$ is not c-weakly
incompressible, it is {\em c-strongly compressible}.

Suppose $K$ is a knot in a closed orientable 3--manifold $M$. We
say that a surface $Q$ is a bridge surface for $K$ if $Q$ is a
Heegaard surface for $M=X \cup_Q Y$ and $K$ intersects each of the
handlebodies $X$ and $Y$  in arcs that are simultaneously parallel
to $Q$. More generally, if the manifold $M$ has boundary, we
require that $K$ intersect each of the compression bodies $X$ and
$Y$ in arcs that are simultaneously parallel into $Q$. The arcs
are called {\em bridges} and the disks of parallelism are called
{\em bridge disks}. A handlebody $H$ intersecting a knot $K$ in a
collection of bridges will be called a $K$-handlebody and will be
denoted by $(H,K)$.

A {\em spine} $\Sigma_H$ of a handlebody $H$ is any graph that $H$
   retracts to.  Removing a neighborhood of a spine
   from a handlebody results in a manifold that is homeomorphic to
   $surface \times I$. More generally the spine $\Sigma_C$ of a compression
   body $C$ is the union of $\bdd_-C$ and a 1--complex $\Gamma$ so
   that $C$ collapses to $\Sigma_C$.

   Suppose $(H,K)$ is a $K$-handlebody and let $\kappa_i$,
   $i=1,\ldots,n$ be the bridges.
   The spine $\Sigma_{(H,K)}$ of
   $(H,K)$ is the union of a
   spine of the handlebody $H$, $\Sigma_H$, together with a collection of
   straight arcs $t_i$, where
   one endpoint of each $t_i$ lies in $\kappa_i$ and the other
   endpoint lies in $\Sigma_H$. Then $H_K -\Sigma_{(H,K)}\cong (\bdd H)_K \times I$.  As in the
   handlebody case, spines of $K$-handlebodies are not unique.

A bridge surface $Q$ for $K$ is called {\em stabilized} if there
is a pair of compressing disks for $Q$, one in $X$ and one in $Y$,
that intersect in exactly one point. The surface $Q$ is called
{\em meridionally stabilized} if there is a compressing disk for
$Q$ in $X$ and a cut-disk for $Q$ in $Y$ (or vice versa) that
intersect in exactly one point. Finally $Q$ is called {\em
perturbed} if there is a pair of bridge disks, $E_X \subset X$ and
$E_Y \subset Y$, such that $E_X \cap E_Y=p$ and $p \in K$ . By
\cite{SchTo071} if $Q$ is stabilized, meridionally stabilized or
perturbed, then there is a bridge surfaces $Q'$ for $K$ such that
$\chi (Q'_K) > \chi (Q_K)$. Furthermore it is easy to see that if
$Q$ is stabilized, meridionally stabilized or perturbed then it is
c-strongly compressible, see \cite{To072}.

Suppose $Q$ is a bridge surface for a knot $K \subset M$.  We say
that $K$ is {\it removable} with respect to $Q$ if it can be
isotoped to lie in $Q$ and there is a meridian disk for one of the
two handlebodies $X$ and $Y$ that intersects $K \subset Q$ in a
single point. In this case, by Lemma 3.3 in \cite{SchTo072}, there
is a Heegaard splitting $\tilde Q$ for $M_K$ such that $\chi
(\tilde Q_K) > \chi (Q_K)$.

We will also need the following result, shown in \cite{To072}.

\begin{thm} \label{thm:essentialexists}
 Suppose $M$ is a closed orientable irreducible 3--manifold containing a
 knot $K$. If $Q$ is a
 c-strongly compressible bridge surface for $K$
  then one of the following is satisfied:

\begin{itemize}

    \item $Q$ is stabilized
    \item $Q$ is meridionally stabilized
    \item $Q$ is perturbed
    \item $K$ is removable
    \item $M$ contains an essential meridional surface $F$ such that $2-\chi(F_K) \leq
2-\chi(Q_K)$.

 \end{itemize}

 \end{thm}

\section{Preliminary results}

Recall that we are considering a Heegaard splitting $A_K \cup P B$
of a knot complement $M_K$ where $M$ is a closed manifold. How
surfaces in a manifold restrict distance of a Heegaard splitting
in various settings has been studied in several papers, for
example see \cite{BS}, \cite{SchTo072} and \cite{ST}. We will take
advantage of some of these results.

\begin{prop} \cite[Proposition 2.5]{Sc2} \label{prop:Proximity2.5}
    Suppose $M$ is an irreducible 3--manifold, $N$ is a compressible
    boundary component of $M$ and $(F, \bdd F)\subset (M, \bdd M)$ is
    a properly embedded essential surface containing no disk
    components and with at least one essential component incident
    to $N$. Let $\mathcal{V}$ be the set of essential curves in $N$
    that bound disks in $M$ and let $f$ be any component of $\bdd F$.
    Then either $d(\mathcal{V}, f) \leq 1-\chi(F)$ or $f$ lies in the
    boundary of $\bdd$-parallel annulus component of $F$.

 \end{prop}

 We will use the above proposition in a very specific situation. The manifold $M$
 will either be the handlebody $B$ or the compression body $A_K$ and $N$ will be $\bdd_+ A_K=\bdd
 B=P$. Then the above proposition says roughly that the distance
 between the boundary curves of an essential surface in a
 compression body and the boundaries of the compressing disks for
 the compression body is bounded above by the Euler characteristic of
 the essential surface.

 The situation becomes considerable more complicated if we allow the
 surface $F$ to have compressing disks. However, if we restrict our
 attention to weakly incompressible surfaces, useful information about
 $d(\mathcal{V}, f)$ can still be obtained, \cite{Sc2}, \cite{To071}. First we recall
 the definition of a tube-spanned recessed collar, \cite{Sc2}. Let $S_0,
 S_1$ be two connected compact subsurfaces in the same component $N$
 of $\bdd M$, with each component $\bdd S_i$, $i=0,1$ essential in
 $\bdd M$ and $S_0 \subset interior(S_1)$. Let $T_i$, $i=0,1$ be the
 properly embedded surface in $M$ obtained by pushing $S_i$, {\it rel} $\bdd
 S_i$ into the interior of $M$, so the region $R$ lying between $T_0$
 and $T_1$ is naturally homeomorphic to $S_1 \times I$. The properly
 embedded surface obtained by tubing $T_0$ and $T_1$ along an
 $I$-fiber of $S_1 \times I$ that is incident to $T_0$ is called a
 {\em tube-spanned recessed collar} in $M$. The properties of these
 surfaces are described in detail in \cite{Sc2}. It turns out that
 tube spanned recessed collars are the only weakly incompressible
 surfaces that don't carry information about distance. More
 precisely:

 \begin{thm}\cite[Theorem 5.4]{Sc2}\label{thm:proximity5.4}
     Suppose $M$ is an irreducible 3--manifold, $N$ is a compressible
    boundary component of $M$ and $(F, \bdd F)\subset (M, \bdd M)$ is
    a bicompressible, weakly incompressible splitting surface with a
    bicompressible component incident to $N$.

    Let $\mathcal{V}$ be the set of essential curves in $N$
    that bound disks in $M$ and let $f$ be any component of $\bdd F \cap
    N$.
    Then either
    \begin{itemize}

    \item $d(\mathcal{V}, f) \leq 1-\chi(F)$ in the curve complex of $N$ or
    \item $f$ lies in the
    boundary of $\bdd$-parallel annulus component of $F$ or
    \item one component of $F$ is a tube spanned recessed collar; all
    other components of $F$ incident to $N$ are incompressible and
    $\bdd$-parallel.

    \end{itemize}
\end{thm}
 Again, we will only be
     considering the situation when $M$ is the handlebody $B$ or the
     compression body $A_K$. Note that in these cases, if the surface $F$
     consists of a single component which is a tube-spanned recessed collar with boundary in $\bdd_+ A_K$
     or in $\bdd B$, then there is
     a spine for $A_K$ or $B$ that is entirely disjoint from $F$ and $F$
     has a compressing disk (a meridional disk for the tube) that lies on
     the same side of $F$ as the spine and is disjoint from the spine.

The above two results tell us that if there is a surface with
certain properties in a compression body, then the boundary curves
are ``not far" from the boundaries of the compressing disks for
the compression body. Thus if $U \cup_R V$ is a Heegaard splitting
for a manifold, $\mathcal{U}$ and $\mcV$ are the collections of
simple closed curves in $R$ that bound disks on sides $U$ and $V$
    respectively and $S$ is a surface that intersects $R$ in a particular way, then we can
hope to obtain a bound on the distance of $R$ by summing $d(\mcU,
R\cap S)$ and $d(R\cap S, \mcV)$. The next Lemma makes this idea
precise.

\begin{lemma} \cite[Lemma 2.6]{ST} \label{lem:essentialdist}
    Suppose $R$ is a Heegaard splitting for a compact manifold $M$, dividing $M$ into compression bodies $U$ and
    $V$. Let $\mathcal{U}$ and $\mcV$ be the collections of simple
    closed curves in $R$ that bound disks on sides $U$ and $V$
    respectively.  Let $S \subset M$ be a properly embedded connected surface
    transverse to $R$ and let $S^U =S \cap U$,
     $S^{V} = S \cap V$. Suppose $S$ satisfies the following
    conditions:
    \begin{itemize}
     \item All curves of $S \cap R$ are essential in $S$ and $R$.

     \item
     There is at least one curve $u \in S \cap R$ such that
     $d(u,\mathcal{U}) \leq 1-\chi(S^U)$ and any curve in $S\cap R$
     for which the inequality does not hold is the boundary of an
     annulus component of $S^U$ that is parallel into $R$.

\item There is at least one curve $v \in S \cap R$ such that
     $d(v,\mathcal{V}) \leq 1-\chi(S^V)$ and any curve in $S\cap R$
     for which the inequality does not hold is the boundary of an
     annulus component of $S^V$ that is parallel into $R$.
    \end{itemize}

Then $d(R) \leq 2-\chi(S)$ .

     \end{lemma}
The following lemma is similar to Lemma 3.6 in \cite{To071} where
the result was proven in the context of a $K$-handlebody. We will
need the result for a compression body so some modifications in
the proof are needed.

\begin{lemma}\label{lem:disjointdisk}

    Let $M$ be a closed orientable irreducible 3--manifold
    containing a knot $K$. Suppose $P$ and $Q$ are Heegaard surface for $M=A \cup_P B=X \cup_Q Y$.
    In addition, suppose that $P$ is also a Heegaard surface for $M_K=A_K \cup_P B$ and $Q$ is a bridge surface for
    $K$. If there is a spine $\Sigma_B$ for $B$ lying in $Y$ and
    $Q_K$ has a c-disk $D$ in $Y_K-P$ that is disjoint from
    $\Sigma_B$ then either
\begin{itemize}
    \item $Q$ is c-strongly compressible, or
    \item $K$ is removable with respect to $Q$, or
    \item $M=S^3$ and $K$ is the unknot.

    \end{itemize}
 \end{lemma}

\begin{proof}
As already noted, $B-\Sigma_B$ has a natural product structure.
Use this structure to push $Q_K$ and $D$ to lie entirely in $A_K$.
As $X \subset A_K$, $Q_K$ always compresses in $X \cap A_K$ and
thus we may assume that $(Q_K,\bdd Q_K) \subset (A_K, \bdd A_K)$
is a c-weakly incompressible surface.

{\bf Case 1:} Suppose first that $D$ is a disk. In this case $Q_K$
is a weakly incompressible surface lying in the compression body
$A_K$. Maximally compress $Q_K$ in $Y_K\cap A_K$ and let $Q'_K$ be
the resulting surface. Note that by the construction, $Q'_K$
separates $P$ and $Q_K$. It is a classical result that maximally
compressing a weakly incompressible surface results in an
incompressible surface, see for example \cite[Lemma 5.5]{Sc2}. As
$Q'_K$ is an incompressible surface in a compression body and
$\bdd Q_K \subset \bdd_-A_K$, each component of $Q'_K$ must be
parallel to $\bdd A_K$, that is, each component is an annulus or a
torus. As $Q'_K$ separates $P$ and $Q_K$ some component $Q^0_K$ of
$Q'_K$ also separates $Q_K$ and $P$.

{\bf Subcase 1A:} Suppose first that $Q^0_K$ is an annulus. Then
the corresponding closed component $Q^0$ of $Q'$ bounds a ball in
$A$ containing $Q$. As $Q'_K$ is parallel to $N(K)$, $K$
intersects this ball in a trivial arc. Now consider $Q^0_K$ as a
surface in the $K$-handlebody $Y$. As handlebodies are
irreducible, $Q^0_K$ must then also bound a ball in $Y$ and by
Lemma 3.2 of \cite{SchTo071}, the knot intersects this ball in a
trivial arc. Thus $M$ is the three sphere and $K$ is a one bridge
knot with respect to the bridge sphere $Q^0$, thus $K$ is the
unknot.

{\bf Subcase 1B:} Suppose then that $Q^0$ is a torus bounding a
solid torus $V$ in $A$ which is a regular neighborhood of $K$. As
$Q_K'$ is obtained from $Q_K$ by compressing along a collection of
disks, the original surface $Q_K$ can be recovered from $Q_K'$ by
tubing along a graph $\Gamma$ with edges dual to these disks.

If $g(Q)=1$, $Q^0$ is the only torus component of $Q'_K$ and all
edges of $\Gamma$ have endpoints on two different components of
$Q'_K$. Pick a meridional disk $F$ for $Q^0$ that intersect $K$
exactly once. Isotope $K$ to lie in $Q^0$ using the parallelism of
$Q^0$ and $Fr(N(K))$, this isotopy is not proper. After the
isotopy of $K$ each edge of $\Gamma$ has at least one endpoint on
a sphere that bounds a ball in $V$. By shrinking $\Gamma$ we can
guarantee that the graph is disjoint from $F$. Thus $F$ is also a
compressing disk for $Q$ that intersects $K$ in a single point
along its boundary, i.e., $K$ is removable with respect to $Q$.

If $g(Q)>1$, consider the two components of $V-Q_K$. One component
is $X_K$ and so it is a handlebody intersecting the knot in
trivial arcs. The other component can be obtained by attaching
1--handles corresponding to the edges of $\Gamma$ to a color of
$Q_K'$. The result is a compression body $Y'$ which intersects $K$
in arcs that are parallel into $\bdd_+ Y'$. Thus $X \cup_Q Y'$ is
a bridge surface for $K$ lying in the manifold $V$. As $g(Q_K)
\geq 2$, by \cite{HS} this implies that $X \cup_ Q Y'$ is
stabilized or perturbed. Thus $Q_K$ is c-strongly compressible in
$M$ as desired.

{\bf Case 2:} Suppose now that $D$ is a cut disk. As $A_K$ is a
compression body, there is a vertical annulus $\Lambda$ in $A$
with one boundary component in $K$ and the other in $P$. As $Q
\cap K \neq \emptyset$, it follows that $Q \cap \Lambda \neq
\emptyset$ and in particular $Q$ intersects $\Lambda$ in at least
one arc. As $Q\cap P=\emptyset$ any such arc must have both of its
endpoints in $K$. Choose $\Lambda$ so that $|\Lambda \cap Q|$ is
minimized, in particular any circle of intersection that is
inessential in $\Lambda$ is essential of $Q$. Now choose $D$ so
that $|D \cap \Lambda|$ is minimal. By a standard innermost circle
argument we may assume that $D \cap \Lambda$ consists only of
circles parallel to the core of $\Lambda$ and arcs, see Figure
\ref{fig:annulusa}.

\begin{figure}
\centering
\includegraphics[width=.8\textwidth]{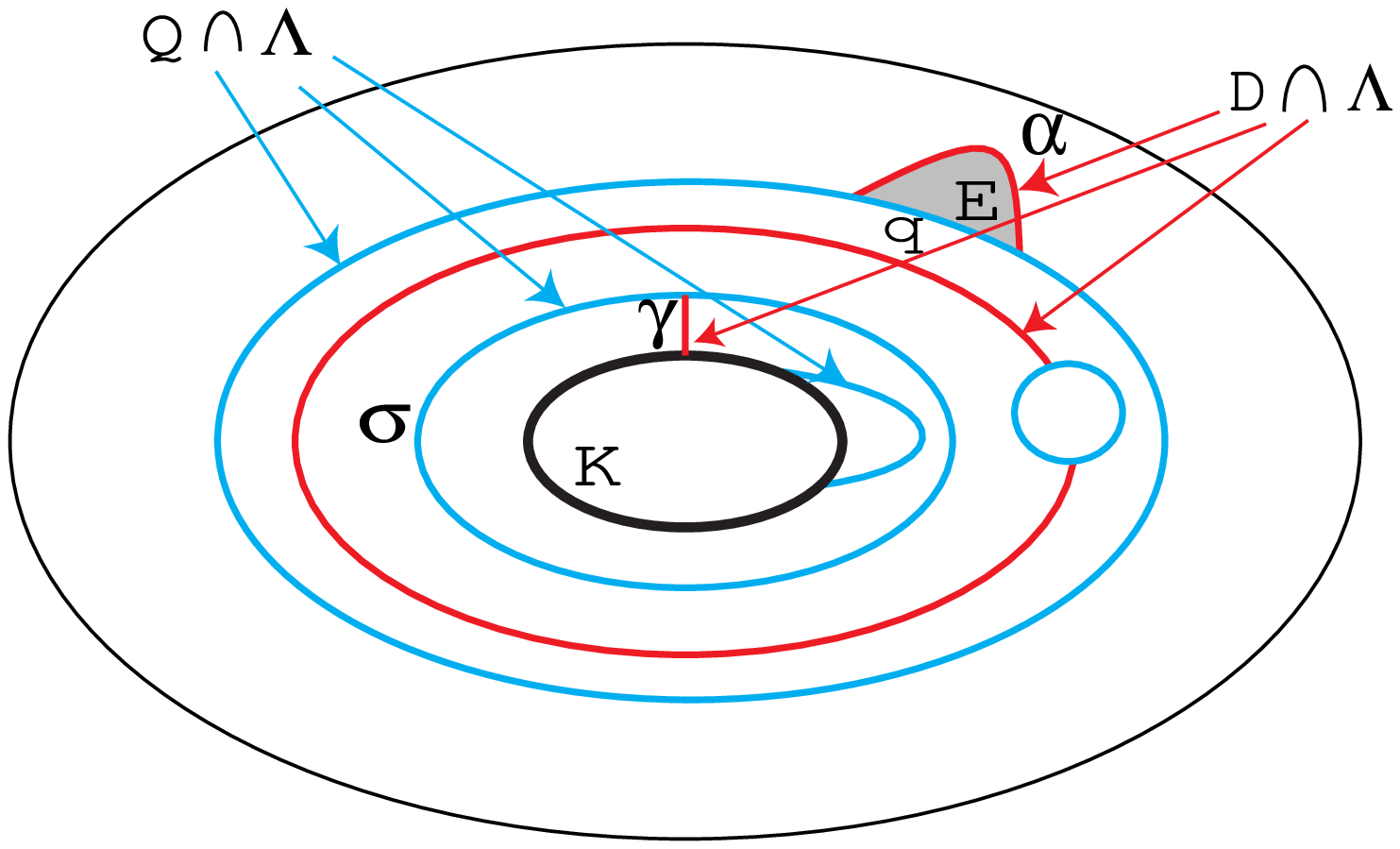}
\caption{} \label{fig:annulusa}
\end{figure}

{\bf Claim 1:} We may assume that every component of $Q \cap
\Lambda$ is either a circle parallel to the core of $\Lambda$ or
is adjacent to an arc of $D \cap \Lambda$.

Suppose $Q \cap \Lambda$ contains either an inessential circle in
$\Lambda$ or an arc with both endpoints of $K$. In either case
assume that this curve is disjoint from $D$. Pick $\delta$ to be
either an innermost such circle or an outermost such arc and let
$G$ be the disk in $\Lambda$ that $\delta$ bounds. By minimality
of $Q \cap \Lambda$ if $\delta$ is a circle, $G$ is a compressing
disk for $Q$. If $\delta$ is an arc, then $G$ is a bridge disk and
so there is a compressing disk for $Q$ contained in the boundary
of a regular neighborhood of $G$. In either case we have found a
compressing disk for $Q$ contained in $A_K$ and disjoint from $P$.
If this disk is contained in $X$, then $Q_K$ is c-strongly
compressible as the disk is disjoint from $D$. If the disk is
contained in $Y$, the result follows by case 1.

\vspace{.3in}

{\bf Claim 2:} There is no arc of $D \cap \Lambda$ that is
parallel in $\Lambda$ to a subarc of $Q \cap \Lambda$.

Suppose there is such an arc and let $\alpha$ be an outermost one,
i.e., $\alpha$ has both endpoints in a component $q$ of $Q$ such
that $q \cup \alpha$ bounds a disk $E$ in $\Lambda$ whose interior
is disjoint from $D$, see Figure \ref{fig:annulusa}. By Claim 1,
the interior of $E$ is also disjoint from $Q$. The boundary of a
regular neighborhood of $D \cup E$ contains at least one c-disk
$F$ for $Q$, see Figure \ref{fig:annulusb}. This c-disk has at
least one fewer intersections with $\Lambda$ than $D$ does. If $F$
is a cut-disk, this will contradict the choice of $D$, if $F$ is a
compressing disk, then Case 1 can be applied. This concludes the
proof of Claim 2.

\vspace{.3in}

\begin{figure}
\centering
\includegraphics[width=.8\textwidth]{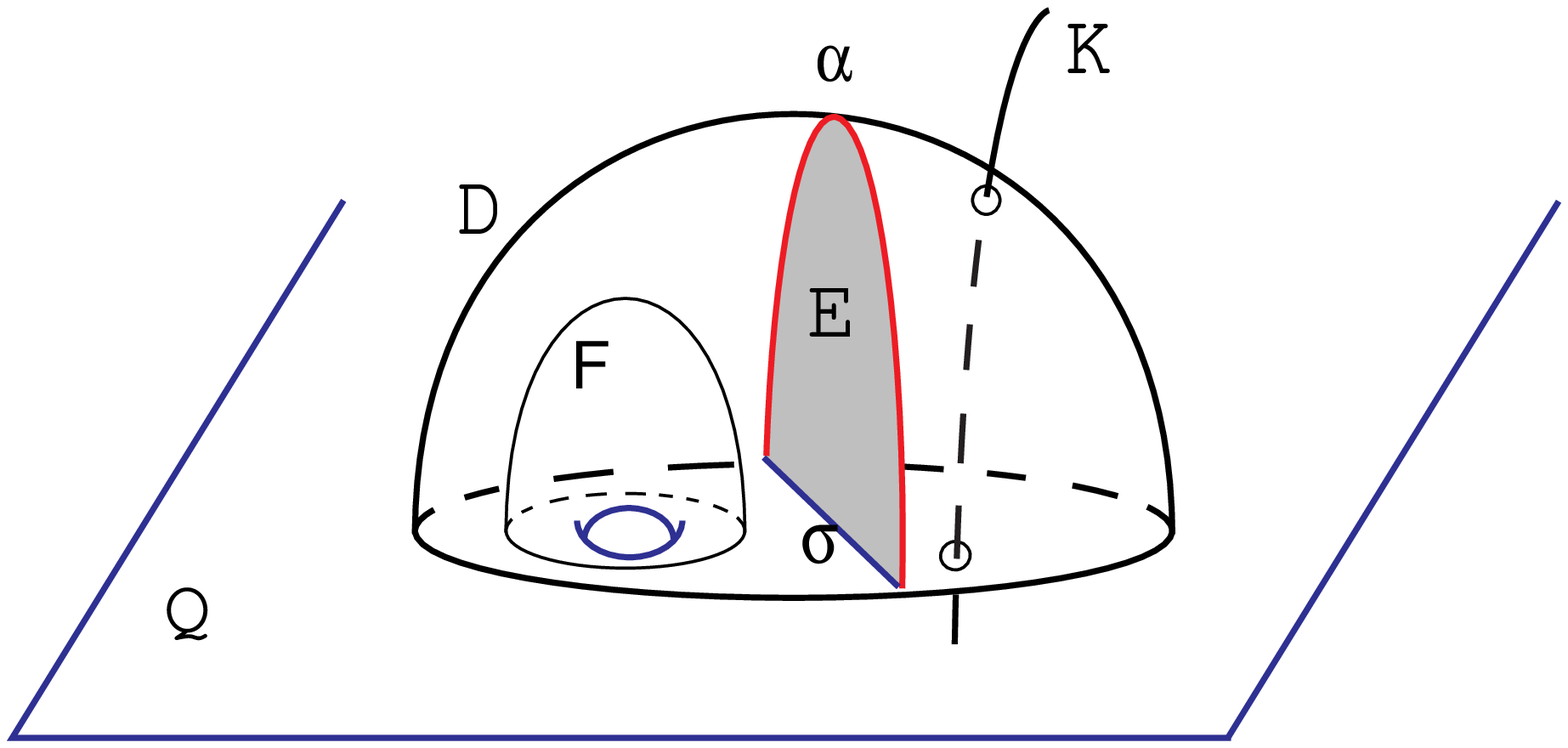}
\caption{} \label{fig:annulusb}
\end{figure}

As $D$ is a cut disk, it intersects the knot exactly once. Let
$\gamma$ be the arc of $D \cap \Lambda$ that has one endpoint in
$K$ and the other in some component $\sigma$ of $Q \cap \Lambda$.

{\bf Subcase A} $\sigma$ is circle parallel to the core of
$\Lambda$, see Figure \ref{fig:annulusa}.

As $Q\cap K \neq \emptyset$, there must be some arc $\delta$ of $Q
\cap \Lambda$ with both endpoints in $K$. This arc is contained in
the disk bounded by $K$, $\sigma$ and two copies of $\gamma$. By
Claim 2, $\delta$ is disjoint from $D$ and so Claim 1 completes
the proof of this case.

{\bf Subcase B} $\sigma$ is a circle that is inessential in
$\Lambda$.

Let $E$ be the disk in $\Lambda$ that $\sigma$ bounds. By Claim 2,
the interior of $E$ is disjoint from $D$ and so by Claim 1 we may
assume that the interior of $E$ is also disjoint from $Q$. Then
$E$ is a compressing disk for $Q_K$ lying in $X$ that intersects
$D$ in exactly one point. By Theorem 2.1 of \cite{To072} $Q_K$ is
c-strongly compressible.

\begin{figure}
\centering
 \includegraphics[width=1\textwidth]{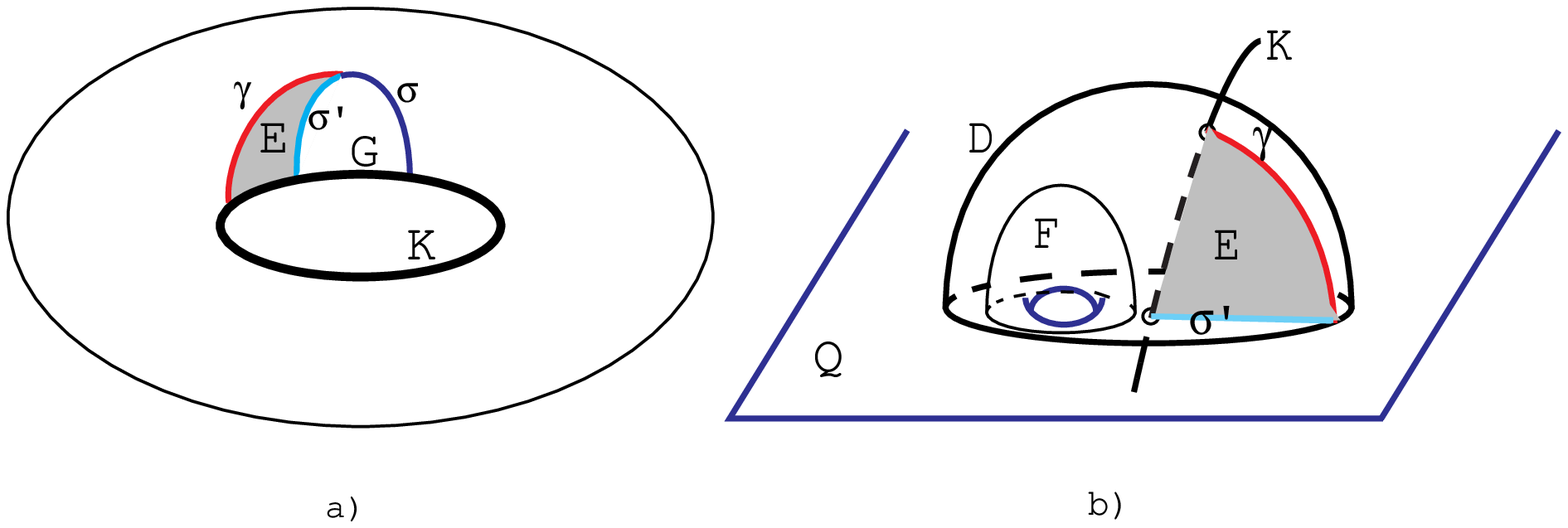}
    \caption{} \label{fig:arc}
   \end{figure}

{\bf Subcase C} $\sigma$ is an arc, see Figure \ref{fig:arc}a.

Let $G$ be the disk that $\sigma$ cuts from $\Lambda$. As before
the boundary of a regular neighborhood of $G$ contains a
compressing disk for $Q$ in $A_K$. If this disk lies in $Y$, Case
1 can be applied so we may assume $G$ lies in $X$. Thus the arc
$\gamma$ lies outside of $G$. Let $E$ be the disk cobounded by
$\gamma$, a subarc $\sigma'$ of $\sigma$ and a subarc of $K$. Then
the boundary of a regular neighborhood of $D \cup E$ contains a
cut disk parallel to $D$ and a disk $F \subset Y$, see Figure
\ref{fig:arc}b. If $\bdd F$ bounds a disk in $Q$, then $\bdd D$
bounds a once punctured disk contradicting the fact that $D$ is a
cut-disk. If $\bdd F$ bounds a punctured disk $D'$ in $Q$, the $F
\cup D'$ is a sphere that intersects the knot once, a
contradiction. Thus $\bdd F$ is essential in $Q_K$ and so
$F\subset Y $ is a compressing disk for $Q_K$ in $A_K$ therefore
we can apply Case 1.

\end{proof}

Finally, the following easy proposition will be used repeatedly in
the proofs to follow.

 \begin{prop} \label{prop:cutimpliescompressing}
      Suppose $P$ is a Heegaard surface for $M_K=A_K \cup_P B$ with genus at least 2 and $D^c$ is a cut disk for $P$ in $A_K$.
      Then there is a compressing disk $D$ for $P$ in $A_K$ such that $d(\bdd D^c, \bdd D)= 1$.

     \end{prop}

\begin{proof}
    Consider cut-compressing $A_K$ along $D^c$, i.e. remove a small open
    neighborhood of $D^c$ from $A_K$. The resulting 3--manifold is a
    $K$-handlebody containing a single bridge - the result of cutting
    the closed loop $K$. As $g(P)\geq 2$, the
    $K$-handlebody has genus at least 1 and thus there is a
    compressing disk $D$ which can be taken to be disjoint from
    the bridge.  This disk is the desired compressing disk for $P$
    in $A_K$.
    \end{proof}

\section{Configurations of $P$ and $Q$}\label{sec:config}

Let $K$ be a nontrivial knot in $M$ and consider how $P$ (a
Heegaard splitting of $M_K=A_K \cup_P B$) and $Q$ (a bridge
surface for $K$ contained in $M=X \cup_Q Y$) intersect in $M$. Let
$Q^A=Q\cap A$, $Q^B=Q\cap B$, $P^X=P\cap X$ and $P^Y=P\cap Y$ and
let $\mcA$ and $\mcB$ be the collections of curves in $P$ that
bound compressing disks in $A$ and $B$ respectively.

\begin{defin} \label{def:removable} Suppose $S$ and $T$ are two
      properly embedded surfaces in a 3--manifold $M$ containing a
      knot $K$ and assume $S$ and $T$ intersect the knot
      transversely. Let $c\in S_K\cap T_K$
     be a simple closed curve bounding possibly punctured disks
      $D^*\subset S_K$ and $E^* \subset T_K$.  If $D^*$ intersects $T_K$ only in curves that are inessential in
      $T_K$ and
      $E^*$ intersects $S_K$ only in curves that are inessential in
      $S_K$ we say that $c$ is {\em removable}.

      \end{defin}

      \begin{rmk} \label{rmk:remove}
\begin{enumerate}

    \item All removable curves can be removed via an isotopy of the
    surfaces that does not affect any essential curves of intersection.

 \item If all curves of intersection between $P$ and $Q_K$ are
 either essential in both surfaces or inessential in both surfaces
 then all inessential curves are removable.
    \end{enumerate}

\end{rmk}

      Isotope $P$ and $Q$ so as to remove all removable curves of
      intersection. We will associate to a position
of $P$ and $Q$ one or more of the following labels.

\begin{itemize}

     \item Label $A$ (resp $B$) if some component of $Q_K \cap
     P$ is the boundary of a compressing disk for $P$ lying in
$A_K$ (resp $B$).

     \item Label $A^c$ if some component of $Q_K \cap
     P$ is the boundary of a cut disk for $P$ lying in $A_K$. (As
     $B \cap K=\emptyset$, no label $B^c$ can occur).

 \item Label $X$ (resp $Y$) if there is a compressing disk for
$Q_K$ lying in $X_K$ (resp $Y_K$) that is disjoint from $P$ and the
configuration does not already have labels $A$, $A^c$ or $B$.

     \item $X^c$ (resp $Y^c$) if there is a cut disk for $Q_K$ lying
     in $X_K$ (resp $Y_K$) that is disjoint from $P$ and the
configuration does not already have labels $A$, $A^c$ or $B$.

     \item $x$ (resp $y$) if some spine $\Sss_{A_K}$ or $\Sss_{B}$ lies
     entirely in $Y_K$ (resp $X_K$) and the configuration does not
already have labels $A$, $A^c$ or $B$.

     \end{itemize}
We will use the superscript $^*$ to denote the possible presence of
superscript $^c$, for example we will use $A^*$ if there is a label
$A, A^c$ or both.

\begin{lemma}\label{lem:nounlabelled}
    If a configuration of $P$ and $Q$ has no labels, then $d(P)\leq
2-\chi(Q_K)$.

\end{lemma}

\begin{proof}
    Suppose a configuration has no labels. First note that this
    implies $P \cap Q_K \neq \emptyset$ as if $P$ is entirely contained
    in $X$ say, then the region would have a label $y$.
    If there is a curve of $P \cap Q_K$ that is inessential in $Q_K$ but
    essential in $P$ an innermost such curve in $Q_K$ would give rise to
    label a $A^*$ or $B$.  If there is a curve of $P \cap Q_K$ that is inessential in $P$ but
    essential in $Q_K$ an innermost such curve in $P$ would give rise to
    label a $X^*$ or $Y^*$. Thus all curves of $P \cap Q_K$ are either
    essential in both surfaces or inessential in both. By Remark
    \ref{rmk:remove}
    we may assume all curves
    of $P \cap Q_K$ are essential in both surfaces.
    Furthermore, as no labels $X^*$ or $Y^*$ are associated to this
    configuration, $Q^A_K$ and $Q^B$ are both c-incompressible and as
    there are no labels $x$ or $y$, they each have at least one
    component that is not parallel to $P$. By Proposition
    \ref{prop:Proximity2.5} both $Q^A_K$ and $Q^B$ satisfy the
    hypothesis of Lemma \ref{lem:essentialdist} from which we deduce
    that $d(P)\leq 2-\chi(Q_K)$.

\end{proof}

\begin{lemma} \label{lem:noAandB} If a configuration of $P$ and $Q$ has labels $A^*$ and $B$,
then $d(P)\leq 2$.
  \end{lemma}

  \begin{proof}
     If a configuration has labels $A$ and $B$, there are curves of
     $P\cap Q_K$, both essential in $P$, that bound compressing disks
     for $P$ in $A_K$ and $B$. As $P$ and $Q_K$ are embedded surfaces,
    these two curves are disjoint and thus $d(P)\leq 1$. If a
    configuration has labels $A^c$ and $B$, using the triangle
    inequality and Proposition \ref{prop:cutimpliescompressing}, we
    deduce that $d(P) \leq 2$.
      \end{proof}

\begin{lemma} \label{lem:noxandY}
If a configuration has labels $x$ and $Y^*$, then at least one of
the following holds:
\begin{itemize}
\item $Q_K$ is c-strongly compressible, or \item $K$ is removable
with respect to $Q$, or \item $d(P)\leq 2-\chi(Q_K)$.
\end{itemize}
\end{lemma}

\begin{proof} Assume $Q_K$ is c-weakly incompressible and $K$ is not removable. We will
show that $d(P)\leq 2-\chi(Q_K)$.
 From the label $x$ we may assume, with no loss of
     generality, that there exists a spine for $A_K$ or $B$ contained in
     $Y_K$.  Note that the spine of $A_K$ contains the frontier of an open neighborhood of $K$ and $K \cap Q
     \neq \emptyset$ thus we conclude that $\Sss_B \subset Y_K$. From
     the label $Y^*$ we know that $Q_K$ has a c-disk in $Y_K - P$,
     call this disk $E$.  By Lemma \ref{lem:disjointdisk}, $E \cap
     \Sss_{B}\neq \emptyset$ so in particular $E \subset B$. Note
     that this also implies $Q_K \cap P\neq \emptyset$.

     As no label $A^*$ or $B$ is present, any curve of intersection that
     is essential in $P$ must also be essential in $Q_K$.
     Suppose there is a curve that is essential in $Q_K$ and inessential
     in $P$. Let $D$ be the disk in $P$ an innermost such curve bounds.
     By Lemma \ref{lem:disjointdisk} the disk $D$ cannot be
     in $Y_K$. If $D$ is in $X_K$, then it is disjoint from $E$ giving a
     c-strong compression for $Q_K$ contrary to our assumption. Thus we
     may assume that all curve of intersection are essential in both $P$
     and $Q_K$.

     Consider first $Q_K^A$.  It is incompressible in $A_K$ because a
     compression into $Y_K$ would violate Lemma \ref{lem:disjointdisk}
     and a compression into $X_K$ would provide a c-weak compression
     of $Q_K$.  If $Q_K^A$ is not essential in $A_K$ then every
     component of $Q_K^A$ is parallel into $P$ so in particular $Q^A\cap
     K =\emptyset$. As it is always the case that $Q^B \cap K =
     \emptyset$ this implies that $Q \cap K =\emptyset$. As $Q$ is a bridge surface for $K$ that is not
     possible.
      We conclude that
     $Q_K^A$ is essential in $A_K$ so by Proposition
     \ref{prop:Proximity2.5} for each component $q$ of $Q_K \cap P$
     that is not the boundary of a $P$-parallel annulus in $A_K$,
     the inequality $d(q, \mathcal{A}) \leq 1 - \chi(Q_K^A)$ holds, i.e., $Q_K^A$
     satisfies the hypotheses of Lemma
     \ref{lem:essentialdist}.

    By Lemma \ref{lem:disjointdisk} $Q^B$
     does not have compressing disks in $Y_K \cap (B - \Sss_B)$ so it
     either has
     no compressing disks in $B - \Sss_B$ at all
     or has a compressing disk lying in $X_K$.   If $Q_B$ is incompressible in $B - \Sss_B
     \cong P \times I$, then either it satisfies the hypothesis of
     Proposition \ref{prop:Proximity2.5} or each component of $Q_B$ is parallel to a subset of
     $P$.  In the latter case the compressing disk $E$ of $Q^B$ in $Y_K- P$ can be
     extended via this parallelism to give a compressing disk for $P$ that is
     disjoint from all $q \in
     Q_K \cap P$.  Hence $d(q, \mathcal{B}) \leq 2\leq 1 -
     \chi(Q^B)$ as long as $Q^B$ is not a collection of
     $P$-parallel annuli.  If that is the case, then $d(\bdd E,
     q_0)$=0 for at least one $q_0 \in (P \cap Q_K)$ so $d(q_0,
     \mathcal{B})\leq 1 \leq 1-\chi(Q^B)$ as desired. Thus in this case
     $Q^B$ satisfies the hypothesis of Lemma \ref{lem:essentialdist} and
     we obtained the desired distance bound.

     If $Q^B$ is bicompressible in $B$, then $Q^B$ is weakly incompressible as every compressing disk for $Q^B$ is also a compressing disk for $Q$.
      Because every compressing disk
     for $Q^B$ in $Y_K$ intersects $\Sigma_B$, $Q^B$ cannot be a
     tube-spanned recessed collar. By Theorem \ref{thm:proximity5.4} $Q^B$ again satisfies the
     hypothesis of Lemma \ref{lem:essentialdist} and we can deduce
     the desired distance bound.
     \end{proof}

     \begin{lemma} \label{lem:noXandY}
    If a configuration has labels $X^*$ and $Y^*$, then at least one of
the following holds:
\begin{itemize}
\item $Q_K$ is c-strongly compressible, or \item $K$ is removable
with respect to $Q$, or \item $d(P)\leq 2-\chi(Q_K)$.
\end{itemize}

     \end{lemma}

\begin{proof} By the labelling we can conclude that $Q_K$ has c-disks in both $X$ and $Y$ that are disjoint from $P$.
Again we will assume that $Q_K$ is c-weakly incompressible and $K$
is not removable and show that the distance bound must hold.
  As before we may assume that all curves of $P\cap Q_K$ are essential
  in both surfaces. A curve that is essential in $P$ but inessential in $Q_K$
  would lead to a label $A^*$ or $B$; a curve that is inessential in
  $P$ but essential in $Q_K$ would lead to c-strong compression of $Q_K$
  as $Q_K$ has c-disks in both $X$ and $Y$ with boundaries disjoint
  from $P$. Thus all curves are either essential or inessential in
  both surfaces and all curves of intersection inessential in
  both surfaces have already been removed.

  The c-disks for $Q_K$ giving rise to the labels $X^*$ and $Y^*$ must
  both be contained in $A_K$ or $B$ by c-weak incompressibility of
  $Q_K$. Suppose $Q_K^A$ is bicompressible and $Q_B$ is
  incompressible, the other case is similar. We may assume that $Q_B$ is essential as otherwise a
  label $x$ or $y$ would be present and Lemma \ref{lem:noxandY} would
  apply. Thus $Q^B$ satisfies the hypothesis of Lemma \ref{lem:essentialdist}.
  On the other hand $Q_K^A$ is a bicompressible c-weakly
  incompressible surface in a 3--manifold. By Theorem \ref{thm:proximity5.4}, $Q_K^A$
  satisfies the hypothesis of Lemma \ref{lem:essentialdist} unless $Q_K^A$ has one
  component that is a tube spanned recessed collar and all other
  components are boundary parallel. It this case $Q_K^A$ would be
  disjoint from a spine of $A_K$ and thus the region would have a
  small label. The result then follows by Lemma \ref{lem:noxandY}

\end{proof}

\begin{lemma}\label{lem:noxandy}
  If a configuration
    has labels $x$ and $y$, then at least one of
the following holds:
\begin{itemize}
\item $Q_K$ is c-strongly compressible, or \item $K$ is removable
with respect to $Q$, or \item $d(P)\leq 2-\chi(Q_K)$.
\end{itemize}
    \end{lemma}

\begin{proof}

Again assume that $Q_K$ is c-weakly incompressible and $K$ is not
removable. We can also assume that the region does not have any
capital labels as otherwise we can apply Lemma \ref{lem:noxandY}.
    We may assume that all curves of intersection are essential in
        both $P$ and $Q_K$ for a curve that is essential in $P$ and
        inessential in $Q_K$ would give rise to a label $A^*$ or $B$ and
        a curve that is essential in $Q_K$ and inessential in $P$ would
        give rise to a label $X^*$ or $Y^*$. Moreover, as labels $X^*$ or $Y^*$
        don't occur, $Q_K^A$ and $Q^B$
        are c-incompressible in $M_K-P$.

        Every spine of $A_K$ intersects $Q$ so both labels $x$ and $y$ must be
        due to a spine of $B$. Supposed there is a spine $\Sss_B\subset Y_K$ and a spine
        $\Sss'_B\subset X_K$. Let $Q_{0}$ be a component of $Q^B$ that lies between
          the two spines.  This implies that $Q_{0}$ is
          parallel into $P$ on both its sides, i.e., that $B \cong Q_{0}
          \times I$. As $g(P)\geq 2$, $Q_0$ is not an annulus.

Let $\aaa$ be an essential arc in $Q_{0}$ with
           endpoints in $P$.
          Then $\aaa \times I \subset Q_{0} \times I \cong B$ is a
          meridian disk $D$ for $B$ that intersects $Q_{0}$ precisely in
          $\aaa$. Consider the frontier of a regular neighborhood of $(Q_0
          \cup D)\cap P$. As $g(P)\geq 2$ the frontier contains at least
          one curve $\sigma$ that is essential in $P$.
          We can conclude
          that for every curve $q \in B\cap Q_K$, $d(\mathcal{B}, q) \leq
          d(\mathcal{B}, \sigma)+d(\sigma,q)\leq 2 \leq 1-\chi(Q^B)$.  Thus $Q^B$ always satisfied
          the hypothesis of Lemma \ref{lem:essentialdist}.

           As we already saw $Q_K^A$ is essential so by Proposition
           \ref{prop:Proximity2.5} $Q_K^{A}$ also satisfies the hypothesis of
           Lemma \ref{lem:essentialdist} and we are done by that lemma.
\end{proof}

Let $\aA$ represent any subset of the labels $A, A^c$, $\bB$ the
label $B$, $\xX$ any subset of the labels $X, X^c, x$ and $\yY$
any subset of the labels $Y, Y^c, y$. Then the lemmas in this
section can be summarized as follows.

\begin{thm} \label{thm:exactlyone}
Suppose $Q_K$ is a c-weakly incompressible bridge surface for a
knot $K$ such that $K$ is not removable with respect to $Q$ and
$P$ is a Heegaard surface for the knot exterior. Let labels $\aA,
\bB, \xX$ and $\yY$ be defined as above. Then either every
configuration of $P$ and $Q_K$ has exactly one of
 the labels $\aA, \bB, \xX$ and $\yY$ associated to it or $d(P)\leq
 2-\chi(Q_K)$.
\end{thm}

\section{Two-Parameter sweep-outs and their graphics}

Let $\Sigma_{(X,K)}$ and $\Sigma_{(Y,K)}$ continue to denote the
spines for the $K$-handlebody $X$ and $Y$. Then there is a map $H:
(Q, Q \cap K) \times I \to (M, K)$ that is a homeomorphism except
over $\Sigma_{(X,K)}\cup \Sigma_{(Y,K)}$ and, near $Q \times \bdd
I$, the map $H$ gives a mapping cylinder structure to a
neighborhood of $\Sigma_{(X,K)} \cup \Sigma_{(Y,K)}$.
 Little is lost and some brevity gained if we restrict $H$ to
 $Q_K \times (I, \bdd I) \to (M_K, \Sigma_{(X,K)}\cup \Sigma_{(Y,K)})$.
The map $H$ is then called a {\em sweep-out} associated to $Q$.
Similarly there is a sweep-out associated to the Heegaard surface
$P$ between the spines $\Sigma_{A_K}$ and $\Sigma_B$.

Consider simultaneous sweep-outs of $P$, between $\Sss_{A_K}$ and
$\Sss_B$ and of $Q_K$ between $\Sss_{(X,K)}$ and $\Sss_{(Y,K)}$.
This two-parameter sweep-out can be described by a square where
each point in the interior of the square represents a position of
$P$ and $Q_K$. Inside the square is a graphic $\Gamma$ which
represents all points where the intersection of the two surfaces
is not generic. At each edge of $\Gamma$ there is a single
tangency between $P$ and $Q_K$ and at each valence 4 vertex there
are two tangencies. Each component of the complement of $\Gamma$
will be called a region and to each region we can associate labels
as defined in Section \ref{sec:config}. Two regions are adjacent
if they share an edge. For the moment we will restrict our
attention to the case when $Q_K$ is a c-weakly incompressible
surface and $K$ is not removable with respect to $Q$, thus by
Theorem \ref{thm:exactlyone} we may assume that either $d(P)\leq
2-\chi(Q_K)$ or every region has exactly one of the labels $\aA,
\bB, \xX$ and $\yY$. In the former case we have achieved the
desired distance bound.

\begin{lemma}\label{lem:noadjacentAandB}
 If two adjacent regions are labelled $\aA$ and $\bB$, then $d(P)\leq 2$.
    \end{lemma}

    \begin{proof}
As the two regions are adjacent, we can transform the
configuration caring the labels $\aA$ into the configuration
caring the label $\bB$ by isotoping $P$ through a single tangency
of $Q_K$. Notice that the curves of intersection of $P$ and $Q_K$
before the isotopy are disjoint from the curves of intersection
after the isotopy. Thus there is a c-disk for $A_K$, $D^*_A$ and a
compressing disk for $B$, $D_B$ such that $D^*_A \cap
D_B=\emptyset$. If $D^*_A$ is actually a compressing disk, that
would imply that $d(P)\leq 1$. If $D^*_A$ is a cut-disk, then we
apply Proposition \ref{prop:cutimpliescompressing} and the
triangle inequality to conclude that $d(P) \leq 2$.

\end{proof}

\begin{lemma} \label{lem:noadjacentXandY}
    Suppose a configuration is changed by passing through a saddle point,
     and the bigon $C$ defining the saddle tangency lies in $X_K \cap
     A_K$ (The case when the bigon lies in the handlebody $B$ is
     similar).  Then

     \begin{itemize}
     \item  No label $x$ or $X^*$ is removed.

     \item  No label $y$ or $Y^*$ is created.

     \item If there is no label $x$ or $X^*$ before
     the move, but one is created after and if there is a label $y$ or
    $Y^*$
     before the move and none after, then $d(P) \leq 2 - \chi(Q_K)$.

     \end{itemize}

     \end{lemma}

    \begin{proof}
    Much of the argument here parallels the argument in the proof of
    Lemma 4.1 in \cite{ST}.

    We first show that no label $x$ or $X^*$ is removed.  If
    there is a c-disk for $X_K$ that lies in $A_K$, a standard
    innermost disk, outermost arc argument on its intersection with $C$
    shows that there is a c-disk for $X_K$ in $A_K$ that is disjoint
    from $C$.  The saddle move has no effect on such a disk (nor,
    clearly,
    on a c-disk for $X_K$ that lies in $B$).  If there is a spine
    of $A_K$ or of $B$ lying entirely in $Y_K$ then that spine, too, is
    unaffected by the saddle move.

    Dually, no label $y$ or $Y^*$ is created: the inverse saddle move,
    restoring the original configuration, is via a bigon that lies in
    $B
    \cap Y_K$.

    To prove the third item position $Q_K$ so that it is exactly tangent
    to $P$ at the
    saddle.  A bicollar of $Q_K$ then has ends that correspond to the
    position of $Q_K$ just before the move and just after.  Let $Q_K^{a}$
    denote $Q_K \cap A_K$ after the move and $Q_K^{b}$ denote $Q_K \cap
    B$ before
    the move.  The bicollar description shows that $Q_K^{a}$ and $Q_K^{b}$
    have disjoint boundaries in $P$.  Moreover the complement of
    $Q_K^{a}
    \cup Q_K^{b}$ in $Q_K$ is a regular neighborhood of the singular
    component
    of $P \cap Q_K$, with Euler characteristic $-1$.  It follows that
    $\chi(Q_K^{a}) + \chi(Q_K^{b}) = \chi(Q_K) + 1$.

    With $Q_K$ positioned as described, tangent to $P$ at the saddle
    point
    but otherwise in general position, consider the closed (non-singular)
    curves of intersection.

    {\bf Claim 1:} It suffices to consider the case in which all
    non-singular curves of intersection are essential in $P$.

    First note that any curve of intersection that is inessential in
    $Q_K$ must be inessential in $P$ as no labels $A^*$ or $B$ are
    present. To prove the claim, suppose a non-singular curve is inessential
    in $P$ and
    consider an innermost one.  Assume first that the disk
    $D$ that it
    bounds in $P$ does not contain the singular curve $s$.  If $\bdd D$ is essential in $Q_K$, then it would
    give
    rise to a label $X$ or a label $Y$ that persists from before the
    move
    until after the move, contradicting the hypothesis.  Thus any
    inessential curve in $P$ that doesn't bound a disk containing the singular curve $s$ is also inessential in $Q_K$ and can be
    removed without affecting the label of the region.

    Suppose then that the disk $D\subset P$ contains the
    singular component $s$. By the above, $s$ is the only
    component of $Q_K \cap P$ in the interior of $D$.
    When the saddle is pushed through, the number of components in $s$
    switches from one, $s_{0}$, to two, $s_{\pm}$, or vice versa.
    All three curves are inessential in $P$ since they lie in $D$.
    The curve $s_0$ and at least one of $s_{\pm}$ bound disks in $P$
    whose interiors are disjoint from $Q_K$. If one of these
    curves was essential in $Q_K$ that would give rise to a label $X$
    or $Y$ that persists through the isotopy. As no such label
    exists, both of these curves are inessential in $Q_K$, ie bound a
    possibly one-punctured disks in $Q_K$. As no
    sphere in $M$ can intersect $K$ exactly once, $s_0$ and one of
    $s_{\pm}$ bound disks in $Q_K$. Because the three curve cobound a pair
    of pants in $Q_K$, all three curves $s_0$ and $s_{\pm}$ are
    inessential in $Q_K$. This means that all three curves are
    removable so passing through this saddle cannot have an effect on
    the labelling.

    {\bf Claim 2:} It suffices to consider the case in which also all
    three curves $s_{0}, s_{\pm}$ are essential in $P$ .

    The case in which all three curves are inessential in $P$ is covered
    in the proof of Claim 1.  If two are inessential in $P$ so is the
    third. Thus the only remaining case is that exactly one of the
    curves $s_{0}, s_{\pm}$ is inessential in $P$ and by Claim 1, the
    disk it bounds in $P$ is disjoint from $Q$. As before the curve
    cannot be essential in $Q_K$ as otherwise it will give rise to a
    label $X$ or $Y$ that persists through the isotopy. Thus the
    curve is inessential in $Q_K$ also (in fact it must bound a disk
    there) so it is removable. If this curve is $s_{\pm}$, passing
    through the saddle can have no effect on the labelling. If the
    removable curve is $s_0$, then the curves $s_{\pm}$ are parallel
    in both surfaces. Passing through the saddle has the
    same effect on the labelling as passing an annulus component of $P^X$
    across
    a parallel annulus component $Q^{0}$ of $Q_K^{A}$.  This move can
    have
    no effect on labels $x$ or $y$.  As there is a label $Y^*$ before
    the move, there is a meridian, possibly punctured disk
    $E^*$ for $Y_K$ that is
    disjoint from $P$. This disk would persist after the move, unless $\bdd E^*$
    is in
    fact the core curve of the annulus $Q^{0}$.  But then the union of
    $E^*$
    and half of $Q^{0}$ would be a possibly punctured meridian disk of
    $A_K$ bounded by a
    component of $\bdd Q^{0} \subset P$.  In other words, there would
    have
    to have been a label $A^*$ before the move, a final contradiction
    establishing Claim 2.

    The above two claims allow us to assume that all curves of
    intersection before and after the move are essential in both
    surfaces. Note that $Q_K^a$ and $Q_K^b$ are c-incompressible (as
    there are no labels $X^*$ or $Y^*$ persisting through the move)
    and have at least one component that is not parallel to $P$ (as
    there are no labels $x$ or $y$ persisting through the move).
    Now apply Proposition \ref{prop:Proximity2.5} to both sides: Let $q_{a}$
    (resp $q_{b}$) be a boundary component of an essential component of
    $Q_K^a$ (resp $Q_K^b$).  Then $$d(P) =
    d(\mathcal{A}, \mathcal{B}) \leq
    d(q_{a}, \mathcal{A}) + d(q_{a}, q_{b}) + d(q_{b}, \mathcal{B}) $$
    $$= \leq 3 -
    \chi(Q_K^a) - \chi(Q_K^b)=2 - \chi(Q_K)$$ as required.

\end{proof}

\begin{lemma} \label{lem:vertexlabels}
    If some vertex of $\Gamma$
    is surrounded by regions labelled with all four labels $\aA, \bB,
    \xX$ and $\yY$, then at least one of
the following holds:
\begin{itemize}
\item $Q_K$ is c-strongly compressible, or \item $K$ is removable
with respect to $Q$, or \item $d(P)\leq 2-\chi(Q_K)$.
\end{itemize}
 \end{lemma}

\begin{proof}
    Suppose there is such a vertex of $\Ggg$ and assume that $Q_K$ is c-weakly incompressible and $K$ is not removable with respect to $Q$. As we have already
    established the desired distance bound if any region has more
    than one label or if two adjacent regions are labelled $\aA$ and
    $\bB$ or $\xX$ and $\yY$ we may assume that going around the
    vertex the regions are labelled in the order $\aA, \xX, \bB$
    and $\yY$. Note then that only two saddle  moves are
    needed to
    move from a configuration labelled $A^*$ to one labelled $B$.  The
    former configuration includes a c-disk for $P$ in $A_K$ and the latter
    a compressing disk for $P$ in $B$. Recall that all curves of
    intersection before a saddle move are disjoint from all curves of
    intersection after the saddle move. Using Proposition
    \ref{prop:cutimpliescompressing} and the triangle inequality, it follows that $d(K,P) \leq 3
    \leq 2 - \chi(Q_K)$, as
    long as at least one of the regions labelled $\xX$ and $\yY$ contains at least one essential curve.

    Suppose all curves of $P\cap Q_K$ in the regions with labels $\xX$ and $\yY$ are inessential.
    Consider the region labelled $\xX$. Crossing the edge in the graphic from this region to the region labelled
    $\aA$ corresponds to attaching a band $b_A$ with both endpoints in an inessential curve
  $c \in P\cap Q$. Note that attaching this band must produce an
  essential curve that gives rise to the label $\aA$, call this curve $c_A$. Similarly crossing
  the edge from the region $\xX$ into the region $B$ corresponds to attaching a band $b_B$
  to give a curve $c_B$. The two bands must be attached to the same
  curve $c$ otherwise
  $c_A$ and $c_B$ would be disjoint curves giving rise to labels $\aA$
  and $\bB$. As we assumed that in the region labelled $\yY$ all
  curves of intersection are inessential,
  attaching both bands simultaneously results in an inessential curve
  $c_{AB}$. But that can only occur of $P$ is a torus which we know
  not to be the case.
\end{proof}

\begin{lemma} \label{lem:edges}

   If a label $\bB$ appears in the regions adjacent to the $\Sigma_{A_K}$ side of $I
    \times I$, or a label $\aA$ appears along the
    $\Sigma_{B}$ side, or a label $\yY$ appears along
    the $\Sigma_{(X,K)}$ side, or a label $\xX$ appears
    along the $\Sigma_{(Y,K)}$ side, then at least one of
the following holds:
\begin{itemize}
\item $Q_K$ is c-strongly compressible, or \item $K$ is removable
with respect to $Q$, or \item $d(P)\leq 2-\chi(Q_K)$.
\end{itemize}
\end{lemma}

\begin{proof}
   Suppose $Q_K$ is a c-weakly incompressible surface and $K$ is not removable with respect to $Q$.
   Consider first a region that is adjacent to the $\Sigma_{A_K}$ side of $I
    \times I$. In such a region $P$ is the boundary of a small
    neighborhood of $\Sigma_{A_K}$ and $Q$ either intersects it in
    meridional circles or doesn't intersect it at all. In the former case any curve of
    intersection $P\cap Q_K$ would lead to a label $\aA$ so if a
    label $\bB$ also appears, then $d(P)\leq 2$ by Lemma
    \ref{lem:noAandB}. In the latter case no label
   $\bB$, or in fact label $\aA$, is possible. Similarly if $P$
    is near $\Sss_B$ and the configuration carries a label $\aA$ then $d(P)\leq 2$.

    Suppose now that $Q_K$ is near $\Sss_{(X,K)}$. As $P$ intersects
    $\Sss_{(X,K)}$ transversely, all but a finite
    number of points of $\Sss_{(X,K)}$ will be disjoint from $P$. Thus a label $\xX$
    necessarily occurs.  If a label $\yY$ also occurs, then $d(P)\leq
    2-\chi(Q_K)$ by Lemma \ref{lem:noXandY}. Symmetrically the
    distance bound holds is some region adjacent to the
    $\Sigma_{(Y,K)}$ boundary of the square is labelled $\xX$

\end{proof}

Finally we will make use of the following combinatorial result.

\begin{thm}\label{thm:sperner}[Sperner's Lemma]
    Suppose a square $I \times I$ contains a graph $\Gamma$ such that
    all vertices of $\Gamma$ in the interior of the square are valence
    4 and all vertices contained in the boundary of the square are
    valence 1.
    Suppose each
    component of $I \times I - \Gamma$ is labelled with exactly one of
    the labels $N, E, S,$ or $W$ in such a way that

    \begin{enumerate}

    \item  no region on the East side of $I \times I$ is labelled $W$,
    no region
    on the West side is labelled $E$, no region on the South side is
    labelled $N$ and no region on the North side is labelled $S$.
    \item no two adjacent regions are labelled $E$ and $W$ nor
 $N$ and $S$.
    \end{enumerate}

    Then some valence 4 vertex is surrounded by regions carrying all 4
    labels.

\end{thm}

\begin{thm} \label{thm:c-weaklyincomp}
    Let $K$ be a knot in a closed irreducible 3--manifold
    $M$ and let $P$ be a Heegaard splitting
    of the knot complement such that $g(P) \geq 2$.
    Suppose $Q$ is a bridge surface for $K$. Then at least one of
the following holds:
\begin{itemize}
\item $Q_K$ is c-strongly compressible, or \item $K$ is removable
with respect to $Q$, or \item $d(P)\leq 2-\chi(Q_K)$.
\end{itemize}
\end{thm}
  \begin{proof}
      Suppose that $Q_K$ is c-weakly incompressible and $K$ is not removable with respect to $Q$.
      Consider a 2-parameter sweep-out of $P$ and $Q_K$ and the
      associated graphic $\Gamma \in I \times I$. Label the regions
      of $I \times I - \Gamma$ with the labels $\aA, \bB, \xX$ and
      $\yY$ as described before. By Theorem
      \ref{thm:exactlyone} we may assume that each region has a
      unique label. By Lemmas \ref{lem:noadjacentAandB} and \ref{lem:noadjacentXandY}
      we may assume that no two adjacent
      regions are labelled $\aA$ and $\bB$ or $\xX$ and $\yY$.
      Finally Lemma \ref{lem:edges} shows that the labels of
      regions adjacent to the boundaries of $I \times I$
   satisfy the conditions that no label $\bB$ appears in the regions adjacent
   to the $\Sigma_{A_K}$ side of $I
    \times I$, no label $\aA$ appears along the
    $\Sigma_{B}$ side, no label $\yY$ appears along
    the $\Sigma_{(X,K)}$ side, and no label $\xX$ appears
    along the $\Sigma_{(Y,K)}$ side of the square. By Theorem
    \ref{thm:sperner} there is a valence 4 vertex of $\Gamma$ that is
    surrounded by regions carrying all four labels. The distance bound
    then follows from Lemma \ref{lem:vertexlabels}.

  \end{proof}

Now we can proof the main result in this paper. We recall the
theorem for the convenience of the reader.

\begin{thm1.4}

Let $K$ be a knot in a closed oriented irreducible 3--manifold
 $M$ and let $P$ be a Heegaard splitting
 of the knot complement of genus at least two. Suppose $Q$ is a bridge surface for $K$. Then either
 \begin{itemize}
 \item $d(P)\leq 2-\chi(Q-K)$, or
 \item $K$ can be isotoped to be disjoint from $Q$ so that after the isotopy $Q$ is
 a Heegaard surface for $M_K$ that is isotopic to a possibly
 stabilized copy of $P$.
\end{itemize}
 \end{thm1.4}
  \begin{proof}
  If $Q_K$ is stabilized, meridionally stabilized or perturbed, as
  described in Section \ref{sec:defin}, there is a bridge surface for $K$,
  $Q'$, such that $\chi(Q'_K) \geq \chi(Q_K)$. By possibly replacing $Q_K$ by $Q'_K$ we may assume
  $Q_K$ is not stabilized, meridionally stabilized or perturbed. If
  $Q_K$ is removable, then again by the results in Section \ref{sec:defin}
  there is a Heegaard surface $\tilde{Q}$ for
  $M_K$ that is isotopic to $Q$ in $M$. By hypothesis $\tilde Q$ is not isotopic to a possibly stabilized copy of $P$ so the result
  follows by Theorem \ref{thm:alternate}. If $Q_K$ is c-weakly incompressible, the result follows by
  Theorem
  \ref{thm:c-weaklyincomp}. The only remaining case is that $Q_K$ is a c-strongly
  compressible bridge surface for $K$ that is not stabilized,
  meridionally stabilized, perturbed or removable. By
  Theorem \ref{thm:essentialexists} in this case there is an essential surface $F$ with $\chi(F_K)
  \geq \chi(Q_K)$. Using Theorem \ref{thm:essentialsurface} we
  deduce that $d(P)\leq 2-\chi(F)\leq 2-\chi(Q_K)$ concluding the
  proof.
\end{proof}

In the case when the manifold is $S^3$ we can also eliminate the
restriction on the genus of $P$ by assuming that $K$ is nontrivial
as that implies that the genus of any Heegaard splitting of
$S^3_K$ is at least 2.

\begin{cor}\label{cor:sphere}

Let $K$ be a non-trivial knot in $S^3$ and let $P$ be a Heegaard
splitting
 of the knot complement. Suppose $Q$ is a bridge surface for $K$. Then either
 \begin{itemize}
 \item $d(P)\leq 2-\chi(Q-K)$, or
 \item $K$ can be isotoped to be disjoint from $Q$ so that after the isotopy $Q$ is
 a Heegaard surface for the knot exterior that is isotopic to a possibly
 stabilized copy of $P$.
\end{itemize}
 \end{cor}

 Notice that if $K$ can be isotoped to be disjoint from $Q$ so that $Q$ is
 a Heegaard surface for $M_K$ that is isotopic to a possibly
 stabilized copy of $P$, then $g(Q)\geq g(P)$. The following
 corollary then follows immediately from Theorem \ref{thm:main} if $g(Q)\geq 1$ and
 from Corollary \ref{cor:sphere}
 if $g(Q)=0$ as the only manifold with a genus 0 Heegaard splitting is $S^3$.

  \begin{cor} \label{cor:main}
  Let $K$ be a knot in a closed irreducible 3--manifold
  $M$ and let $P$ be a Heegaard splitting
  of the knot complement. Suppose $Q$ is a bridge surface for $K$ such
  that $g(Q) < g(P)$. Then $d(P)\leq 2-\chi(Q-K)$.
  \end{cor}

\end{document}